\documentclass[a4paper,11pt]{article}
\usepackage{amsmath, amsfonts, amscd, amssymb, amsthm, enumerate, xypic}

\newcommand{\Mat}{\operatorname{M}}
\newcommand{\Mats}{\operatorname{S}}
\newcommand{\Mata}{\operatorname{A}}

\newcommand{\tr}{\operatorname{tr}}
\newcommand{\GL}{\operatorname{GL}}
\newcommand{\Ker}{\operatorname{Ker}}

\renewcommand{\setminus}{\smallsetminus}

\def\F{\mathbb{F}}

\def\lcro{\mathopen{[\![}}
\def\rcro{\mathclose{]\!]}}

\theoremstyle{definition}

\theoremstyle{plain}
\newtheorem{theo}{Theorem}

\newtheorem{step}{Step}

\theoremstyle{plain}

\theoremstyle{remark}

\title{From trivial spectrum subspaces to spaces of diagonalizable real matrices}
\author{Cl\'ement de Seguins Pazzis\footnote{Universit\'e de Versailles Saint-Quentin-en-Yvelines, Laboratoire de Math\'ematiques
de Versailles, 45 avenue des \'Etats-Unis, 78035 Versailles cedex, France}
\footnote{e-mail address: dsp.prof@gmail.com}}

\begin{document}


\thispagestyle{plain}

\maketitle

\begin{center}
\emph{Dedicated to Bernard Rand\'e, on the occasion of his retirement.}
\end{center}

\begin{abstract}
A recent generalization of Gerstenhaber's theorem on spaces of nilpotent matrices \cite{dSPlargeaffine} is shown
to yield a new proof of the classification of linear subspaces of diagonalizable real matrices with the maximal dimension \cite{RandedSP}.
\end{abstract}



\section{Introduction}
Let $\F$ be an arbitrary field and $n$ be a positive integer.
We denote respectively by $\Mat_n(\F)$, $\Mats_n(\F)$ and $\Mata_n(\F)$ the spaces of all $n$ by $n$ square matrices, symmetric matrices, and alternating
matrices with entries in $\F$. We denote by $\GL_n(\F)$ the group of all invertible matrices of $\Mat_n(\F)$.

A \textbf{trivial spectrum subspace} of $\Mat_n(\F)$
is a linear subspace in which no matrix has a non-zero eigenvalue in the field $\F$.
Such subspaces have attracted much focus in the recent past, in part due to their connection with affine subspaces of non-singular
matrices, and also because they are a generalization of Gerstenhaber's linear subspaces of nilpotent matrices \cite{Gerstenhaber}.
A subset $V$ of $\Mat_n(\F)$ is called \textbf{irreducible} when no non-trivial linear subspace of $\F^n$
is stable under all the elements of $V$. The main theorem on trivial spectrum subspaces then reads as follows
(see \cite{dSPlargeaffine} for the original proof, and \cite{dSPQJM} for a simplified proof over fields with large cardinality).

\begin{theo}[de Seguins Pazzis (2013)]\label{TStheo}
Let $V$ be an irreducible trivial spectrum subspace of $\Mat_n(\F)$.
Then:
\begin{enumerate}[(a)]
\item $\dim V \leq \dbinom{n}{2}$.
\item If $\dim V=\dbinom{n}{2}$ and $|\F|>2$ then
$V=P\Mata_n(\F)$ for some matrix $P$ of $\GL_n(\F)$ which is non-isotropic, i.e.\ $\forall X \in \F^n \setminus \{0\}, \; X^TPX \neq 0$.
\end{enumerate}
\end{theo}

\vskip 3mm
In this short note, we shall uncover a surprising connection between trivial spectrum subspaces and large linear subspaces of diagonalizable matrices.
The main theorem on the latter \cite{RandedSP}, which is especially relevant to the case of real matrices, reads as follows:

\begin{theo}[Rand\'e, de Seguins Pazzis (2011)]\label{Diagtheo}
Let $V$ be a diagonalizable subspace of $\Mat_n(\F)$, i.e.\ a linear subspace in which every matrix is diagonalizable.
Then :
\begin{enumerate}[(a)]
\item $\dim V \leq \dbinom{n+1}{2}$.
\item If $\dim V=\dbinom{n+1}{2}$ then $V$ is similar to $\Mats_n(\F)$, i.e.\ $V=P \Mats_n(\F) P^{-1}$ for some $P \in \GL_n(\F)$.
\end{enumerate}
\end{theo}

Statement (a) is readily obtained by noting that $V$ is linearly disjoint from the linear subspace of all
strictly upper-triangular matrices of $\Mat_n(\F)$.

Over an arbitrary field, not all symmetric matrices are diagonalizable (see \cite{Waterhouse}),
and hence the maximal dimension for a diagonalizable subspace of $\Mat_n(\F)$ might be less than $\dbinom{n+1}{2}$.
In the extreme case when $\F$ is algebraically closed with characteristic $0$, the celebrated Motzkin-Taussky theorem \cite{MoTau2}
shows that the maximal dimension is $n$.

We shall prove that, provided that $\F$ has more than $2$ elements, point (b) in Theorem \ref{Diagtheo}
is a rather simple consequence of Theorem \ref{TStheo}.

\section{Proof of Theorem \ref{Diagtheo}}

For $(i,j)\in \lcro 1,n\rcro^2$, we denote by $E_{i,j}$ the matrix unit of $\Mat_n(\F)$
in which every entry equals zero except the one at the $(i,j)$-spot, which equals $1$.

\vskip 3mm
The key idea is to use orthogonality with respect to the non-degenerate symmetric bilinear form
$$(A,B) \in \Mat_n(\F)^2 \mapsto \tr(AB).$$
Note that if we take a linear subspace $V$ of $\Mat_n(\F)$ and an invertible matrix $P \in \GL_n(\F)$,
the orthogonal complement of $PVP^{-1}$ is $P V^\bot P^{-1}$.

Now, let $V$ be a diagonalizable subspace of $\Mat_n(\F)$ with $\dim V=\dbinom{n+1}{2}$.

The proof has four basic steps.

\begin{step}
The subspace $V$ contains the identity matrix $I_n$.
\end{step}

\begin{proof}
We note that $V':=V+\F I_n$ is a diagonalizable subspace that includes $V$. Statement (a) from Theorem \ref{Diagtheo}
yields $\dim V' \leq \dbinom{n+1}{2}=\dim V$, and hence $V'=V$. It follows that $I_n \in V$.
\end{proof}

In the next two steps, we prove that $V^\bot$ is an irreducible trivial spectrum subspace of $\Mat_n(\F)$.

\begin{step}
The space $V^\bot$ is irreducible.
\end{step}

\begin{proof}
Assume that the contrary holds. Then, $n \geq 2$ and by replacing $V$ with a well-chosen similar subspace,
we see that no generality is lost in assuming that we have found an integer $p \in \lcro 1,n-1\rcro$ such that
every matrix of $V^\bot$ reads
$$\begin{bmatrix}
[?]_{p \times p} & [?]_{p \times (n-p)} \\
[0]_{(n-p) \times p} & [?]_{(n-p) \times (n-p)}
\end{bmatrix}.$$
Hence, $V=(V^\bot)^\bot$ contains $E_{1,n}$, which is not diagonalizable, contradicting our assumptions.
\end{proof}

\begin{step}\label{TSclaim}
The space $V^\bot$ is a trivial spectrum one.
\end{step}

\begin{proof}
Assume on the contrary that $V^\bot$ contains a matrix with a non-zero eigenvalue in $\F$. Replacing $V$
with a similar subspace, we can assume that $V^\bot$ contains a matrix of the form
$$H=\begin{bmatrix}
\lambda & [?]_{1 \times (n-1)} \\
[0]_{(n-1) \times 1} & [?]_{(n-1) \times (n-1)}
\end{bmatrix} \quad \text{for some $\lambda \in \F \setminus \{0\}$.}$$
We have linear mappings $a : V \rightarrow \F$, $C : V \rightarrow \F^{n-1}$, $R : V \rightarrow \Mat_{1,n-1}(\F)$ and $K : V \rightarrow \Mat_{n-1}(\F)$
such that every matrix $M$ of $V$ reads
$$M=\begin{bmatrix}
a(M) & R(M) \\
C(M) & K(M)
\end{bmatrix}.$$
Set $W:=\bigl\{M \in V : C(M)=0\bigr\}$.
By the rank theorem,
$$\dim V = \dim C(V)+\dim W\leq (n-1)+\dim W.$$
Let $M \in W$ be such that $K(M)=0$. Since $M$ is orthogonal to $H$, we get $\lambda a(M)=0$ and hence $a(M)=0$.
It follows that $M$ is both diagonalizable and strictly upper-triangular, whence $M=0$. Therefore, $\dim K(W)=\dim W$.

Next, let $M \in W$. Denoting by $e_1$ the first vector of the standard basis of $\F^n$, we see that $X \mapsto MX$ stabilizes
$\F e_1$ and the induced endomorphism on the quotient space $\F^n /\F e_1$ is represented by $K(M)$. Since $M$ is diagonalizable,
so is $K(M)$. Hence, point (a) of Theorem \ref{Diagtheo} applies to $K(W)$ and yields
$$\dim V \leq (n-1)+\dim K(W) \leq (n-1)+\dbinom{n}{2}=\dbinom{n+1}{2}-1,$$
which contradicts our assumptions.
\end{proof}

Finally, we note that $\dim V^\bot=n^2-\dim V=\dbinom{n}{2}$.
Assume now that $|\F| >2$. Then, Theorem \ref{TStheo} yields a non-isotropic matrix $P \in \GL_n(\F)$ such that $V^\bot=P \Mata_n(\F)$,
whence
$$V=\bigl(P \Mata_n(\F)\bigr)^\bot=\Mats_n(\F)\, P^{-1}.$$
Since $I_n \in V$, the matrix $P$ is symmetric.
Note that, for all $Q \in \GL_n(\F)$,
$$\Mats_n(\F)\bigl(Q PQ^T\bigr)^{-1}=QQ^{-1} \Mats_n(\F) (Q^{-1})^T P^{-1} Q^{-1}=Q \Mats_n(\F) P^{-1} Q^{-1}=QV Q^{-1}.$$
We draw two consequences from that remark:
\begin{enumerate}[(1)]
\item Replacing $V$ with a similar subspace amounts to replacing $P$ with a congruent matrix.
\item If $P$ is congruent to a scalar multiple of $I_n$ then $V$ is similar to $\Mats_n(\F)$.
\end{enumerate}
Thus, our final step, which will complete the proof, follows:

\begin{step}
The matrix $P$ is congruent to a scalar multiple of $I_n$.
\end{step}

\begin{proof}
The result is obvious if $n=1$, so we assume that $n \geq 2$.
Since $P$ is symmetric and non-isotropic, it is congruent to a diagonal matrix
(even if $\F$ has characteristic $2$).
By point (1) above, no generality is lost in assuming that $P$ is actually diagonal, with non-zero diagonal entries $d_1,\dots,d_n$.

Fix $i \in \lcro 2,n\rcro$. Then, $(E_{1,i}+E_{i,1}) P^{-1}=d_i^{-1} E_{1,i}+d_1^{-1} E_{i,1}$ must be diagonalizable.
Yet, its characteristic polynomial equals $t^{n-2}(t^2-d_i^{-1} d_1^{-1})$ and hence,
for some $\lambda \in \F \setminus \{0\}$, we successively get $d_i^{-1} d_1^{-1}=\lambda^2$ and
$d_i=d_1(d_1^{-1} \lambda^{-1})^2$ with $d_1^{-1} \lambda^{-1} \neq 0$.
We conclude that $P$ is congruent to $d_1 I_n$.
\end{proof}

A final note on the case of the field with $2$ elements, which we have left out:
For subspaces of $\Mat_2(\F_2)$, the above proof still works because the result of
Theorem \ref{TStheo} is known to hold in that case. Yet, over $\F_2$ the symmetric matrix $\begin{bmatrix}
1 & 1 \\
1 & 1
\end{bmatrix}$ is non-diagonalizable because it is nilpotent and non-zero.
Hence, $\dim V<3$ for every diagonalizable subspace of $\Mat_2(\F_2)$.
From there, an induction shows that $\dim V<\dbinom{n+1}{2}$ whenever $V$ is a diagonalizable subspace of $\Mat_n(\F_2)$
with $n \geq 2$ (with the notation from the proof of Step \ref{TSclaim}, one applies the induction hypothesis to $K(W)$
and one notes that $W \cap \Ker K \cap \Ker a=\{0\}$).

\end{document}